# Breakdown of analyticity for $\bar{\partial}_b$ and Szegö kernels

by

Joe KAMIMOTO

Department of Mathematical Sciences University of Tokyo,
3-8-1, Komaba, Meguro, Tokyo, 153 Japan.
*E-mail adress* : `kamimoto@ms.u-tokyo.ac.jp`

**Abstract**

The CR manifold $M_m = \{\text{Im} z_2 = [\text{Re} z_1]^{2m}\} (m = 2, 3, \ldots)$ is the counterexample, which has been given by M. Christ and D. Geller, to analytic hypoellipticity of $\bar{\partial}_b$ and real analyticity of the Szegö kernel. In order to give a direct interpretation for the breakdown of real analyticity of the Szegö kernel, we give a Borel summation type representation of the Szegö kernel in terms of simple singular solutions of the equation $\bar{\partial}_b u = 0$.

# 0  Introduction

In [10], M. Christ and D. Geller gave the following remarkable counterexample to analytic hypoellipticity of $\bar{\partial}_b$ for real analytic CR manifolds of finite type (in the sense of Kohn [26] or D'Angelo [12]):

**THEOREM 0.1**  *On the three-dimensional CR manifold $M_m := \{\mathrm{Im}\, z_2 = [\mathrm{Re}\, z_1]^{2m}\}$ ($m = 2, 3, \ldots$), $\bar{\partial}_b$ fails to be analytic hypoelliptic (in the modified sense of [10]).*

Moreover analytic hypoellipticity of $\bar{\partial}_b$ is closely connected with real analyticity of the Szegö kernel off the diagonal. By considering the Szegö kernel as a singular solution of the equation $\bar{\partial}_b u = 0$, Christ and Geller obtained Theorem 0.1 as a corollary to the following theorem.

**THEOREM 0.2**  *The Szegö kernel of $M_m$ ($m=2, 3, \ldots$) fails to be real analytic off the diagonal.*

The proof of Theorem 0.2 by Christ and Geller [10] is based on certain formulas of A. Nagel [30]. Though their proof is logically clear, it seems difficult to understand the singularity of the Szegö kernel of $M_m$ directly, since their proof was established by contradiction. On the other hand, M. Christ ([6],[7]) directly constructed singular solutions of $\bar{\partial}_b u = 0$ and proved Theorem 0.1. In this paper, we give an integral representation of the Szegö kernel of $M_m$ in terms of the singular solutions of Christ. Since the singular solutions of Christ are substantially simpler, our integral representation makes it easy to understand the singularity of the Szegö kernel of $M_m$. Moreover we give the direct proof of Theorem 0.2. We also give a similar representation of the Bergman kernel of the domain $\{\mathrm{Im}\, z_2 > [\mathrm{Re}\, z_1]^{2m}\} \subset \mathbf{C}^2$ ($m = 2, 3, \ldots$) on the boundary.

We remark that our subject in this paper has also been treated in the paper of M. Christ ([7],§7), and that our result can be considered as an improvement of Proposition 7.2 in [7].

In general we consider the hypersurface

$$M_P := \{\mathrm{Im}\, z_2 = P(z_1)\} \subset \mathbf{C}^2,$$

where $P : \mathbf{C} \to \mathbf{R}$ is a real analytic function. We assume that $\triangle P$ is nonnegative and does not vanish identically. Such a surface is pseudoconvex and of finite type. A nonvanishing, antiholomorphic, tangent vector field is $\partial/\partial\bar{z}_1 - 2i(\partial P/\partial\bar{z}_1)\partial/\partial\bar{z}_2$. As coordinates for the surface we use $\mathbf{C} \times \mathbf{R} \ni (z = x+iy, t) \mapsto (z, t+iP(z))$; the vector field pulls back to $\bar{\partial}_b = \partial/\partial\bar{z} - i(\partial P/\partial\bar{z})\partial/\partial t$. Let $\bar{\partial}_b^*$ denote the formal adjoint of $\bar{\partial}_b$ with respect to the Lebesgue measure on $\mathbf{C} \times \mathbf{R}$. We say that $\bar{\partial}_b$ is *analytic hypoelliptic on $M_P$ (in the modified sense of* [10]) if whenever $\bar{\partial}_b u$ is real analytic in an open set $U$ and $u = \bar{\partial}_b^* v$ for some $v \in L^2$ in $U$, $u$ is real analytic in $U$. In the usual sense, $\bar{\partial}_b$ is not even $C^\infty$ hypoelliptic, but it is known that if $\bar{\partial}_b \bar{\partial}_b^* u \in C^\infty$, then $\bar{\partial}_b^* u \in C^\infty$ ([27]).

Let $S((z,t);(w,s))$ be the Szegö kernel of $M_P$; that is, the distribution kernel associated to the operator defined by the orthogonal projection of $L^2(\mathbf{C} \times \mathbf{R})$, with respect to the Lebesgue measure, onto the kernel of $\bar{\partial}_b$. It is known that the Szegö kernel is $C^\infty$ off the diagonal ([31],[32]).

There are many interesting studies on real analyticity for $\bar{\partial}_b$ and Szegö kernels. For certain strictly pseudoconvex CR manifolds $M_P$ (i.e. $\triangle P > 0$), $\bar{\partial}_b$ is analytic hypoelliptic ([17]) and the Szegö kernel is real analytic off the diagonal ([22]). In the weakly pseudoconvex case, there are many important results (see the references in [10]), but the precise condition for them to be real analyticity is still unknown.

In investigating the properties of Szegö kernels, it is an important problem to find a *good* expression for them. For some classes of CR manifolds, the Szegö kernels (or the Bergman kernels for some domains) are explicitly computed in closed form ([19],[11],[13],[16], [21]). But the Szegö kernels for almost all CR manifolds cannot be written in closed form. We consider the following two types of the integral representation for the Szegö kernels. They are very useful to analyse the singularities of the Szegö kernel (e.g. [3],[10],[20],[16],[24]).

The representations of first type were obtained in [18],[29],[28],[30],[21]. For example, the Szegö kernel of $M_P$ ($P$ satisfies certain conditions) can be represented as follows (see [21]):

$$S((z,t);(w,s)) = c\int_0^\infty \int_{-\infty}^\infty \frac{\tau \exp(\tau(\eta(z+\bar{w}) - P(z) - P(w) - i(s-t)))}{\displaystyle\int_{-\infty}^\infty \exp(2\tau(r\eta - P(r)))dr} d\eta d\tau.$$

The representations of this type can be obtained by using the generalized Paley-Wiener theorem. Those of second type are represented as the Borel (or Mittag-Leffler) summation of some countably many functions. For example, Bonami and Lohoué [3] gave the representation for the Szegö kernel of $\{\sum_{j=1}^{n} |z_j|^{2m_j} = 1\}$ ($m_j \in \mathbf{N}$):

$$S(z,w) = c \int_0^\infty e^{-p} \left[ \prod_{j=1}^n \sum_{\nu_j=0}^\infty \frac{(z_j \bar{w}_j p^{\frac{1}{m_j}})^{\nu_j}}{\Gamma(\frac{\nu_j}{m_j} + \frac{1}{m_j})} \right] p^{\sum_{j=1}^n \frac{1}{m_j} - 1} dp.$$

In this paper, we give a representation of the Borel summation type for the Szegö kernel of $M_m$ by using the representation of first type, which is obtained by Nagel ([30]). Since the Szegö kernel is expressed by the superposition of certain simple singular solutions of $\bar{\partial}_b u = 0$ due to M. Christ, the structure of the singularity can be understood directly. The singularity of the Szegö kernel is almost equal to that of Christ's singular solution involving the first eigenfunction of the certain ordinary differential operator (see §6).

In this paper, we analyse the counterexample of Christ and Geller directly by using the classical asymptotic analysis for ordinary differential equations with irregular singular points. In particular, some techniques to obtain the asymptotic expansion of the functions admitting an integral representation or a Taylor series expansion are very useful for our computation ([14],[36]). Explicitly the properties of the entire function:

$$\varphi(x) = \int_{-\infty}^\infty e^{-2(w^{2m} - xw)} dw$$

play an important role in the breakdown of real analyticity, so we must obtain the detailed informations of $\varphi$. The properties of the function $\varphi$ are studied in [25] in detail. The function $\varphi$ satisfies the ordinary differential equation of higher order:

$$\frac{d^{2m-1}}{d^{2m-1}x} y - \frac{2^{2m-2}}{m} xy = 0.$$

The above equation is a particular form of so-called *Turrittin's equation*, which has been studied by many authors (see the references in [25]).

The plan of this paper is as follows. We state our results and outline of the proofs in Section 1. In Section 2, we recall the direct construction of the singular solutions of $\bar{\partial}_b u = 0$ in [6],[7], which are used in our representations. In Section 3, we establish our theorems. In Sections 4,5, we give the proofs of propositions and lemmas respectively, which are necessary for the proofs of our theorems. In Section 6, we show the failure of real analyticity of the Szegö kernel directly by using our representation. In Section 7, we study the boundary behaviors of the Szegö and Bergman kernels of the domain $\{\mathrm{Im} z_1 > [\mathrm{Re} z_1]^{2m}\}$ $(m = 2, 3, \ldots)$ on the diagonal.

In this paper, we use $c$ or $C$ ($C(X_1, X_2, \cdots)$) for various constants (depending on $X_1, X_2, \cdots$) without further comment.

# 1 Statement of main results

For $M_m = \{\mathrm{Im} z_2 = [\mathrm{Re} z_1]^{2m}\}$ $(m = 2, 3, \ldots)$, M. Christ in [6],[7] constructed the singular solutions of the equation $\bar{\partial}_b u = 0$ ($u = \bar{\partial}_b^* v, v \in L^2$) by applying the partial Fourier transformation and by solving a certain simple ordinary differential equation (see §2).

Let $\varphi(x)$ be the function defined by

$$\varphi(x) = \int_{-\infty}^{\infty} e^{-2(w^{2m} - xw)} dw.$$

Then it is known that $\varphi$ has infinitely many simple zeros and that all of them exist on the imaginary axis ([33]). We denote them by $\pm i a_j$ ($j \in \mathbf{N}$), where the $a_j$'s are positive and arranged in the increasing order. (More detailed information of $\varphi$ is explained in Section 3.)

Christ's solutions are of the following form:

$$(1.1) \qquad S_j^v(z, t) = \int_0^\infty e^{it\tau} e^{-x^{2m}\tau} e^{\sigma(y) i a_j z \tau^{\frac{1}{2m}}} \tau^v d\tau \quad (y \neq 0, j \in \mathbf{N}, v \geq 0),$$

where $z = x + iy$ and $\sigma(y)$ is the sign of $y$. It is easy to check that the $S_j^v$'s are not real analytic on $\{(0 + iy, 0); y \in \mathbf{R}\}$. Besides this, the $S_j^v$'s, off the set $\{y = 0\}$, belong to $s$th order Gevrey class $G^s$ for all $s \geq 2m$, but no better, where $G^s := \{f; \exists C > 0 \text{ s.t. } |\partial^\alpha f| \leq C^{|\alpha|} \Gamma(s|\alpha|) \ \forall \alpha\}$.

For $M_m$ ($m = 2, 3, \ldots$), define the distribution $K(z, t) = S((z, t); (0, 0))$. Then $K$ is a $C^\infty$ function away from $(0, 0)$ ([31],[32]). We give a representation of $K$ in terms of the singular solutions $\{S_j^v\}_{j \in \mathbf{N}}$.

THEOREM 1.1 *If* $|\arg z \pm \frac{\pi}{2}| < \frac{1}{2m-1} \frac{\pi}{2}$, *then*

$$(1.2) \qquad K(z, t) = c^S \int_0^\infty e^{-p} H(z, t; p) dp,$$

*where*

(1.3) $$H(z,t;p) = \sum_{j=1}^{\infty} c_j S_j^{\frac{1}{m}}(z,t) p^{f(j)},$$

*for some sequence* $f(j) = j + j_0 + O(j^{-1})(> 0)$ *as* $j \to \infty$ *and* $c^S$, $c_j$*'s and* $j_0 \in \mathbf{Z}$ *are constants. Moreover there exists a positive constant* $C^S(z)$ *depending on z such that,*

(1.4) $$|H(z,t;p)| \leq C^S(z)|p|^{-\frac{1}{4}}.$$

Let $B((z_1,z_2);(w_1,w_2))$ be the Bergman kernel of the domain $D=\{\text{Im} z_2 > P(z_1)\} \subset \mathbf{C}^2$; that is, the distribution kernel for the orthogonal projection of $L^2(D)$ onto the subspace of holomorphic functions. Then $B$ extends to a $C^\infty$ function on $\overline{D} \times \overline{D}$ minus the diagonal ([31],[32]). When $D = D_m := \{\text{Im} z_2 > [\text{Re} z_1]^{2m}\}$ $(m = 2,3,\ldots)$, we obtain a similar representation of $K^B(z,t) := B((z,t+ix^{2m});(0,0))$.

THEOREM **1.2** *If* $|\arg z \pm \frac{\pi}{2}| < \frac{1}{2m-1}\frac{\pi}{2}$, *then*

(1.5) $$K^B(z,t) = c^B \int_0^\infty e^{-p} H^B(z,t;p) dp,$$

*with*

(1.6) $$H^B(z,t;p) = \sum_{j=1}^{\infty} c_j S_j^{1+\frac{1}{m}}(z,t) p^{f(j)},$$

*where* $c_j$*'s and* $f(j)$ *are as in Theorem 1.1 and* $c^B$ *is a constant. Moreover there exists a positive constant* $C^B(z)$ *depending on z such that,*

(1.7) $$\left|H^B(z,t;p)\right| \leq C^B(z)|p|^{-\frac{1}{4}}.$$

*Remarks.* 1) The constants $c_j$'s in (1.3), (1.6) are given as follows:

$$c_j = \frac{1}{\varphi'(ia_j)} \frac{1}{\Gamma(f(j)+1)} \quad j \in \mathbf{N}.$$

2) We believe that $j_0$ and $O(j^{-1})$ in the above theorems are removable, that is, $f(j)$ can be replaced by $j$. We do not know whether the order $-\frac{1}{4}$ in (1.4),(1.7) is optimal.

3) If we change the order of the sum and the integral in (1.2),(1.5) formally, we obtain the formal sum of the form $\sum_{j=1}^{\infty} d_j S_j^v(z,t)$, where $d_j$'s are constants. However the formal sum is not convergent in the usual sense. We shall show this fact in Subsection 5.6.

**Outline of the proofs.** We explain the idea of the proof of Theorem 1.1 roughly. The proof of Theorem 1.2 is given in the same fashion.

Our computation starts from the following formula which has been obtained by A. Nagel [30]:

$$(1.8) \quad K(z,t) = c \int_0^\infty e^{it\tau} e^{-x^{2m}\tau} \left[ \int_{-\infty}^\infty e^{z\tau^{\frac{1}{2m}}v} \frac{1}{\varphi(v)} dv \right] \tau^{\frac{1}{m}} d\tau.$$

As mentioned above, the function $\varphi$ has countably many zeros. Therefore by shifting the integral contour formally, applying the residue formula and changing the order of the sum and the integral, we can obtain the formal sum $\sum_{j=1}^{\infty} d_j S_j^{\frac{1}{m}}(z,t)$, where $d_j$'s are constants and $S_j^{\frac{1}{m}}$'s are as in (1.1). But this sum is not convergent in the usual sense (see §5.6). In order to justify this formal computation, we use the idea of the exact WKB method (see e.g. [35],[1]).

Introducing the positive large parameter $q$, we consider the integral

$$(1.9) \quad K_q(z,t) = c \int_0^\infty e^{it\tau} e^{-x^{2m}\tau} \left[ \int_\Gamma e^{z\tau^{\frac{1}{2m}}v} \frac{q^{-F(v)-1}}{\varphi(v)} dv \right] \tau^{\frac{1}{m}} d\tau,$$

where $F$ and $\Gamma$ are an appropriate function and contour, respectively. Note that when $q = 1$, this representation reduces to (1.8). The same procedure as above yields the convergent sum $\sum_{j=1}^{\infty} d_j S_j^{\frac{1}{m}}(z,t) q^{-f(j)-1}$, where $f(j)$ is as in the theorem. Next, we apply the Borel and Laplace transformations to the above convergent sum with respect to $q$ in this order, then we can obtain the following representation:

$$c \int_0^\infty e^{-qp} \left[ \sum_{j=1}^\infty \frac{d_j}{\Gamma(f(j)+1)} S_j^{\frac{1}{m}}(z,t) p^{f(j)} \right] dp,$$

which is equal to the above convergent sum. Now we regard $q$ as a complex variable. Since the sum in the integral has the estimate (1.4), the above integral can be analytically continued to $q = 1$. In fact, the sum in the integral has the estimate (1.4) in the theorem. Hence we obtain the integral representation (1.2) in the theorem.

The estimate (1.4) of the sum $H$ can be obtained by using the idea of E. M. Wright [36]. The method of Wright is useful to obtain the asymptotic expansion of entire functions which are expressed by Taylor series.

## 2 Construction of singular solutions

In this section we recall M. Christ's construction [6],[7] of singular solutions to $\bar{\partial}_b u = 0$.

We show that there exist the function $F$ such that $\bar{\partial}_b^* F$ is not real analytic and $\bar{\partial}_b \bar{\partial}_b^* F \equiv 0$. We identify $M_m$ with $\mathbf{C} \times \mathbf{R}$ as in the Introduction. Let $\bar{\partial}_b^*$ be the formal adjoint of $\bar{\partial}_b$ with respect to the Lebesgue measure on $\mathbf{C} \times \mathbf{R}$.

In the case of $M_m = \{\mathrm{Im} z_2 = [\mathrm{Re} z_1]^{2m}\}$ ($m \in \mathbf{N}$), $\bar{\partial}_b = X + iY$ and $\bar{\partial}_b^* = -X + iY$ where

$$X = \frac{\partial}{\partial x}, \quad Y = \frac{\partial}{\partial y} - 2mx^{2m-1}\frac{\partial}{\partial t}.$$

Applying the partial Fourier transformation in the $y$ and $t$ variables, we seek solutions to $\bar{\partial}_b \bar{\partial}_b^* u \equiv 0$ of the form $u(x,y,t) = e^{it\tau} e^{i\eta y} f(x)$. Then $\bar{\partial}_b \bar{\partial}_b^* u = 0$ reduces to the ordinary differential equation:

$$(2.1) \quad \left[-\frac{d}{dx} + (\eta - 2m\tau x^{2m-1})\right] \circ \left[\frac{d}{dx} + (\eta - 2m\tau x^{2m-1})\right] f(x) = 0.$$

We suppose that $\tau > 0$ and change the variables by $u = \tau^{\frac{1}{2m}} x$, $\eta = \tau^{\frac{1}{2m}} \xi$, with $\xi \in \mathbf{C}$. If we set $f(x) = g(\tau^{\frac{1}{2m}} x)$, then (2.1) reduces to

$$L_\xi g(u) := \left[-\frac{d}{du} + (\xi - 2mu^{2m-1})\right] \circ \left[\frac{d}{du} + (\xi - 2mu^{2m-1})\right] g(u) = 0$$

Since $L_\xi$ factors as a product of two first-order operators, one finds an explicit solution of $L_\xi g_\xi = 0$:

$$g_\xi(u) = e^{-\xi u + u^{2m}} \int_{-\infty}^{u} e^{-2(s^{2m} - \xi s)} ds.$$

Here we fix $\xi \in \mathbf{C}$ and define formally

$$F_\xi(z,t) = \int_0^\infty e^{it\tau} e^{i\xi y \tau^{\frac{1}{2m}}} g_\xi(\tau^{\frac{1}{2m}} x) d\tau.$$

Then $F_\xi$ satisfies $\bar{\partial}_b \bar{\partial}_b^* F_\xi \equiv 0$ formally and

$$\bar{\partial}_b^* F_\xi(z,t) = \int_0^\infty e^{it\tau} e^{-x^{2m}\tau} e^{z\xi \tau^{\frac{1}{2m}}} d\tau.$$

If the imaginary part $\sigma$ of $\xi$ is positive (resp. negative) and if $g_\xi$ is a bounded function on $\mathbf{R}$, then the above two integrals converge absolutely for $y > 0$ ( resp. $y < 0$). Moreover

$$\left| \frac{\partial^k}{\partial t^k} \bar{\partial}_b^* F_\xi(0+iy, 0) \right| = \left| \int_0^\infty \tau^k e^{-\tau^{\frac{1}{2m}} \sigma y} d\tau \right|$$

$$= 2m(\sigma y)^{-2mk-2m} \Gamma(2mk+2m).$$

Thus $\bar{\partial}_b^* F_\xi$ would not be real analytic and moreover $\bar{\partial}_b^* F_\xi$ would belong to $s$-th order Gevrey class $G^s$ for all $s \geq 2m$, but no better, where $G^s = \{f \, ; \, \exists C > 0 \, s.t. \, |\partial^\alpha f| \leq C^{|\alpha|} \Gamma(s|\alpha|) \, \forall \alpha\}$.

From the above, if there exists $\xi \in \mathbf{C}$, with $\mathrm{Im}\,\xi \neq 0$, such that $g_\xi$ is bounded on $\mathbf{R}$, we obtain singular solutions of $\bar{\partial}_b u = 0$. The following lemma gives the condition for $g_\xi$ to be bounded.

LEMMA 2.1 ([5]) *The function $g_\xi$ is bounded on $\mathbf{R}$ if and only if $\xi \in \mathbf{C}$ satisfies*

$$\varphi(\xi) = \int_{-\infty}^\infty e^{-2(w^{2m} - \xi w)} dw = 0.$$

Here we remark that the function $\varphi$ appears in the integral formula of the Szegö kernel of $M_m$ (1.8), which has been computed by A. Nagel [30].

As was mentioned in Section 2, the function $\varphi$ has infinitely many zeros and all of them exist on the imaginary axis for $m \in \{2, 3, \ldots\}$. Hence we obtain non-real analytic solutions, in the range of $\bar{\partial}_b^*$, to $\bar{\partial}_b u = 0$ for $m \in \{2, 3, \ldots\}$. However when $m = 1$, $\varphi$ is a Gaussian function and it has no zeros, so we can not construct the singular solutions.

# 3 Proofs of Theorems 1.1 and 1.2

In this section we give the proofs of Theorems 1.1 and 1.2. We suppose that $m$ is an integer with $m \geq 2$.

## 3.1 Proof of Theorem 1.1

A. Nagel [30] has computed the Szegö kernel of $M_m$ explicitly. We set $K(z,t) = S((z,t);(0,0))$, then we have

$$(3.1) \qquad K(z,t) = c \int_0^\infty e^{it\tau} e^{-x^{2m}\tau} P(z\tau^{\frac{1}{2m}}) \tau^{\frac{1}{m}} d\tau,$$

with

$$(3.2) \qquad P(u) = \int_{-\infty}^\infty e^{uv} \frac{1}{\varphi(v)} dv,$$

and

$$\varphi(v) = \int_{-\infty}^\infty e^{-2(w^{2m} - vw)} dw.$$

From the above formula, in order to investigate the properties of $K(z,t)$, the analysis of the function $\varphi$ seems to be important. Actually the following properties of $\varphi$ in the two lemmas below are very important in our computations. From now on we regard the function $\varphi$ as an entire function on the complex plane. Note that $\varphi$ is an even function and $\varphi(0) > 0$.

The first lemma is concerned with the asymptotic behavior of $\varphi$ at infinity. By $f(x) \sim g(x)$ as $x \to \infty$ on $\alpha \leq \arg x \leq \beta$ ($\alpha \leq \beta$), we mean that $f(x)/g(x)$ converges 1 as $|x| \to \infty$ uniformly on $\alpha \leq \arg x \leq \beta$.

LEMMA 3.1 Let $A(x)$ be defined by $A(x) = c_0 x^{\frac{1-m}{2m-1}} \exp\{c_1 x^{\frac{2m}{2m-1}}\}$, where $c_0 = \pi^{\frac{1}{2}}(2m-1)^{-\frac{1}{2}}(2m)^{-\frac{1}{4m-2}}$ and $c_1 = 2[(\frac{1}{2m})^{\frac{1}{2m-1}} - (\frac{1}{2m})^{\frac{2m}{2m-1}}]$. Then

$$\varphi(x) \sim A(x) + A(xe^{-\pi i}), \qquad \text{as } x \to \infty \text{ on } -\frac{\pi}{2} + \varepsilon < \arg x < \frac{3\pi}{2} - \varepsilon,$$

where $\varepsilon$ is an arbitrary constant with $0 < \varepsilon < \frac{\pi}{4}$.

From the above lemma, we have

$$\varphi(x) \sim A(x), \qquad \text{as } x \to \infty \text{ on } -\frac{\pi}{2} + \varepsilon \leq \arg x \leq \frac{\pi}{2} - \varepsilon,$$

and

$$\varphi(x) \sim A(xe^{-\pi i}), \quad \text{as } x \to \infty \text{ on } \frac{\pi}{2} + \varepsilon \leq \arg x \leq \frac{3}{2}\pi - \varepsilon.$$

The second lemma is concerned with the distribution of the zeros of $\varphi$. Since the exponential order of $\varphi$ is $\frac{2m}{2m-1}$, $\varphi$ has infinitely many zeros. Moreover, it is known in [33] that all zeros of $\varphi$ exist on the imaginary axis. The asymptotic distribution of the zeros of $\varphi$ is as follows:

LEMMA **3.2** *All but finitely many of zeros of $\varphi$ are simple. Let $\{\pm i a_j; 0 < a_j < a_{j+1}, j \in \mathbf{N}\}$ be the set of simple zeros of $\varphi$. Then we have*

$$j = c_2 a_j^{\frac{2m}{2m-1}} - \frac{1}{4} - j_0 + O(j^{-1}), \text{ as } j \to \infty,$$

*where $c_2 = \frac{2}{\pi}[(\frac{1}{2m})^{\frac{1}{2m-1}} - (\frac{1}{2m})^{\frac{2m}{2m-1}}] \cdot \cos \frac{1}{2m-1}\frac{\pi}{2}$ and $j_0 \in \mathbf{Z}$ is a constant depending only on $m$.*

We denote the set of non-simple zeros of $\varphi$ by $\{\pm i b_j; 0 < b_j < b_{j+1}, j = 1, 2, \ldots, J\}$, where $J \geq 0$ is an integer depending only on $m$. We expect that the above set is empty (i.e. $J = 0$).

For the proof of the two lemmas above, we refer to the forthcoming paper ([25]).

For the computation below, we prepare the integral contours $\Gamma_\pm$, $L_\pm$ in the $v$-plane. Set $r_0 = \frac{1}{2}\min\{a_1, b_1\}$. The contour $\Gamma_+$ is defined in the following way: $\Gamma_+$ follows the half-line $v = re^{i\frac{m+1}{m}\frac{\pi}{2}}$ from $r = \infty$ to $r = r_0$, the circle $v = r_0 e^{+i\theta}$ from $\theta = \frac{m+1}{m}\frac{\pi}{2}$ to $\theta = \frac{m-1}{m}\frac{\pi}{2}$ and the half-line $v = re^{\frac{m-1}{m}\frac{\pi}{2}}$ from $r = r_0$ to $r = \infty$. The contour $L_+$ is defined in the following way: $L_+$ consists of two parts $L_+^1$, $L_+^2$, which follow the imaginary axis. $L_+^1$ follows it from $v = +i\infty$ to $v = ir_0$ and $L_+^2$ from $v = ir_0$ to $v = +i\infty$, but $L_+^1$ passes left-hand side of each points $v = ia_j, ib_j$, while $L_+^2$ right-hand side of them. Reflecting the contours $\Gamma_+$, $L_+$ with respect to the real axis, we obtain $\Gamma_-$, $L_-$, respectively. (See Figure 1.)

Now we define $K_q(z, t)$ as follows :

(3.3) $$K_q(z, t) = c \int_0^\infty e^{it\tau} e^{-x^{2m}\tau} P_q(z\tau^{\frac{1}{2m}}) \tau^{\frac{1}{m}} d\tau,$$

Figure 1: Integral contours $\Gamma_\pm$, $L_\pm$.

with

(3.4) $$P_q(u) = \int_{\Gamma_{\sigma(y)}} e^{uv} \frac{q^{-F_{\sigma(y)}(v)-1}}{\varphi(v)} dv,$$

(3.5) $$F_{\pm}(v) = c_2 e^{\mp \frac{2m}{2m-1} \frac{\pi}{2}} v^{\frac{2m}{2m-1}} - \frac{1}{4} = c_2 \left[\pm \frac{1}{i} v\right]^{\frac{2m}{2m-1}} - \frac{1}{4},$$

where $\sigma(y)$ is the sign of $y$ and $q$ is a complex parameter belonging to a region on which the integral in (3.4) makes sense. Such a region is given by the following lemma.

LEMMA **3.3** *There exist positive constants $\varepsilon_0 < 1$ and $\alpha_0 < \frac{\pi}{2}$ such that $P_q(u)$ is a holomorphic function of $q$ in the domain $V := \{q \in \mathbf{C}; |q| > 1 - \varepsilon_0 \text{ and } |\arg q| < \alpha_0\}$.*

From now on we suppose that $q$ belongs to $V$. We remark that $q = 1$ is contained in $V$ and $K_1(z,t) = K(z,t)$. In fact we can deform the integral contour $\mathbf{R}$ into $\Gamma_+$ or $\Gamma_-$ in (3.4) by Lemma 3.1.

First we give the proof of the theorem in the case where $z$ is in the sector $|\arg z - \frac{\pi}{2}| < \frac{1}{2m-1} \frac{\pi}{2}$. By Lemma 3.1, there exists a positive number $q_0$ such that if $|\arg v - \frac{\pi}{2}| < \frac{\pi}{2m}$, then

$$\left| e^{uv} \frac{q_0^{-F_+(v)-1}}{\varphi(v)} \right| \leq C \exp\left\{-c|v|^{\frac{2m}{2m-1}} \cos \frac{\pi}{2m-1}\right\},$$

where $C, c \ (> 0)$ are constants independent of $v$. Thus we have

$$P_{q_0}(u) = \int_{\Gamma_+} e^{uv} \frac{q_0^{-F_+(v)-1}}{\varphi(v)} dv = \int_{L_+} e^{uv} \frac{q_0^{-F_+(v)-1}}{\varphi(v)} dv,$$

by Cauchy's theorem. Moreover by the residue formula and Lemma 3.2, we have

(3.6) $$P_{q_0}(u) = 2\pi i \sum_{j=1}^{\infty} \frac{1}{\varphi'(a_j i)} e^{ia_j u} q_0^{-f(j)-1},$$

for some sequence $f(j)$ such that

(3.7) $$f(j) = j + j_0 + O(j^{-1}) \text{ as } j \to \infty,$$

where $j_0 \in \mathbf{Z}$ is as in Lemma 3.2. Substituting (3.6) into (3.3) and changing the order of the sum and the integral, we have

$$(3.8) \qquad K_{q_0}(z,t) = c \sum_{j=1}^{\infty} \frac{1}{\varphi'(a_j i)} S_j^{\frac{1}{m}}(z,t) q_0^{-f(j)-1},$$

with

$$(3.9) \qquad S_j^v(z,t) = \int_0^{\infty} e^{it\tau} e^{-x^{2m}\tau} e^{\sigma(y) i a_j z \tau^{\frac{1}{2m}}} \tau^v d\tau,$$

where $a_j$'s are as in Lemma 3.2. Note that $S_j^v$, $j \in \mathbf{N}$ and $v > 0$, is a slight generalization of the singular solution $\bar{\partial}_b^* F$ constructed in Section 2.

Now we define by $H(z,t;p)$ the Borel transformation of (3.8) with respect to $q$, where

$$H(z,t;p) = \sum_{j=1}^{\infty} \frac{1}{\Gamma(f(j)+1)} \frac{1}{\varphi'(a_j i)} S_j^{\frac{1}{m}}(z,t) p^{f(j)}.$$

Then we have
$$K_{q_0}(z,t) = c \int_0^{\infty} e^{-q_0 p} H(z,t;p) dp.$$

In the same fashion, we have

$$(3.10) \qquad K_q(z,t) = c \int_0^{\infty} e^{-qp} H(z,t;p) dp,$$

for $q \geq q_0$. In fact if $q \geq q_0$, then

$$\left| e^{uv} \frac{q^{-F(v)-1}}{\varphi(v)} \right| \leq \left| e^{uv} \frac{q_0^{-F(v)-1}}{\varphi(v)} \right|,$$

in $|\arg v - \frac{\pi}{2}| < \frac{\pi}{2m}$. Moreover, since the right-hand side of (3.10) extends analytically with respect to $q$ to the region $\{\operatorname{Re} q > q_0\} \cap V$, (3.10) is satisfied there. If we admit Proposition 3.1 below, then (3.10) is satisfied on the region $V$ in the same fashion.

PROPOSITION 3.1 *If* $|\arg z \pm \frac{\pi}{2}| < \frac{1}{2m-1} \frac{\pi}{2}$, *then there is a positive constant* $C^S(z)$ *depending on $z$ such that*

$$|H(z,t;p)| \leq C^S(z) |p|^{-\frac{1}{4}}.$$

In particular, (3.10) is satisfied when $q = 1$. Hence we have

$$K(z,t) = K_1(z,t) = c \int_0^\infty e^{-p} H(z,t;p) dp.$$

On the other hand, if $z$ is in the sector $|\arg z + \frac{\pi}{2}| < \frac{1}{2m-1}\frac{\pi}{2}$, then we can also obtain (3.4) by replacing $\Gamma_+$ with $\Gamma_-$, $F_+(v)$ with $F_-(v)$ and $+ia_j$ with $-ia_j$ in the above argument.

This completes the proof of Theorem 1.1. $\square$

## 3.2 Proof of Theorem 1.2

The relation between the Szegö kernel of $M_m$ and the Bergman kernel $B((z_1, z_2); (w_1, w_2))$ of the domain $D_m = \{\text{Im} z_2 > [\text{Re} z_1]^{2m}\} \subset \mathbf{C}^2$ is obtained in [32], §7. We define $K^B(z,t) = B((z, t+ix^{2m}); (0,0))$. This relation and the exact computation of the Szegö kernel of $M_m$ by Nagel ([30]) yield

$$(3.11) \qquad K^B(z,t) = c \int_0^\infty e^{it\tau} e^{-x^{2m}\tau} P(z\tau^{\frac{1}{2m}}) \tau^{1+\frac{1}{m}} d\tau,$$

where $P$ is as in (3.2). The difference between (3.1) and (3.11) does not give any essential influence on the argument in Subsection 3.1. Therefore if we admit Proposition 3.2 below, then we can obtain in the same fashion as in Subsection 3.1

$$K^B(z,t) = c^B \int_0^\infty e^{-p} H^B(z,t;p) dp,$$

with

$$H^B(z,t;p) = \sum_{j=1}^\infty \frac{1}{\varphi'(a_j i)} \frac{1}{\Gamma(f(j)+1)} S_j^{1+\frac{1}{m}}(z,t) p^{f(j)},$$

where $f(j)$ is as in (3.7) and $S_j^{1+\frac{1}{m}}$ is as in (3.9).

PROPOSITION **3.2** *If* $|\arg z \pm \frac{\pi}{2}| < \frac{1}{2m-1}\frac{\pi}{2}$, *then there is a positive constant* $C^B(z)$ *depending on $z$ such that*

$$|H^B(z,t;p)| \leq C^B(z)|p|^{-\frac{1}{4}}.$$

$\square$

# 4 Proof of Propositions 3.1 and 3.2

In this section we give the proofs of Propositions 3.1 and 3.2.

Now we define $H^v(z,t;p)$ by

(4.1) $$H^v(z,t;p) = \sum_{j=1}^{\infty} \frac{1}{\Gamma(f(j)+1)} \frac{1}{\varphi'(\sigma(y)ia_j)} S_j^v(z,t) p^{f(j)},$$

where $v \geq 0$ and $S_j^v(z,t)$ are as in (3.9). Note that $H = H^{\frac{1}{m}}$ and $H^B = H^{1+\frac{1}{m}}$. In order to prove the two propositions, it is sufficient to show the following:

PROPOSITION 4.1 *If* $|\arg z \pm \frac{\pi}{2}| < \frac{1}{2m-1}\frac{\pi}{2}$, *then there is a positive constant* $C^v(z)$ *depending on* $z$ *and* $v$ *such that*

$$|H^v(z,t;p)| \leq C^v(z)|p|^{-\frac{1}{4}}.$$

*Proof of Proposition 4.1.* First we consider the case where $z$ is in the sector $|\arg z - \frac{\pi}{2}| < \frac{1}{2m-1}\frac{\pi}{2}$.

By Cauchy's theorem and the residue formula,

(4.2) $$H^v(z,t;p) = \int_{L_+} \frac{1}{\Gamma(F_+(\xi)+1)} \frac{1}{\varphi(\xi)} S^v(\xi;z,t) p^{F_+(\xi)} d\xi,$$

with

$$S^v(\xi;z,t) = \int_0^{\infty} e^{it\tau} e^{-x^{2m}\tau} e^{\xi z \tau^{\frac{1}{2m}}} \tau^v d\tau,$$

where the function $F_+(\xi)$ and the integral contour $L_+$ is as in Section 3. Define the domains $D_\xi^\pm, D_\zeta$ in $\mathbf{C}$ by

$$D_\xi^\pm = \left\{\xi \in \mathbf{C}; \left|\arg \xi \mp \frac{\pi}{2}\right| < \frac{2m-2}{2m-1}\frac{\pi}{2}\right\} \quad \text{and} \quad D_\zeta = \left\{\zeta \in \mathbf{C}; \operatorname{Re}\zeta > -\frac{1}{4}\right\}$$

respectively. Then $F_\pm : D_\xi^\pm \to D_\zeta$ are biholomorphic functions. Let $G_\pm : D_\zeta \to D_\xi^\pm$ be the inverse functions of $F_\pm$. Here $G_\pm(\zeta) = \pm i[\frac{1}{c_2}\left(\zeta + \frac{1}{4}\right)]^{\frac{2m-1}{2m}}$.

Let $A_j$ be the values of $F_+(ia_j) = F_-(-ia_j) \in \mathbf{R}$ for $j \in \mathbf{N}$, $B_j$ the values of $F_+(ib_j) = F_-(-ib_j)$ for $j = 1, 2, \cdots, J$. Let $\mathcal{N}$ be the domain $\{\zeta; |\zeta - A_j| < \delta$ or $|\zeta - B_j| < \delta\}$, where $\delta > 0$ is a sufficiently small number.

We define the integral contours $L, M_1, M_2, N$ on $\overline{D_\zeta}$ in the following way: $L$ consists of two parts $L^1, L^2$, which follow the real axis. $L^1$ follows it from $\infty$ to $F(ir_0)$ and $L^2$ from $F(ir_0)$ to $\infty$, but $L^1$ passes below each points $A_j, B_j$, while $L^2$ passes above such points. $N$ follows the line $\text{Re}\zeta = -\frac{1}{4}$ upwards. $M_1$ follows the circle $\zeta = Re^{i\theta}$ from $\theta = -\theta_1$ to $\theta = \theta_1$, where $\cos\theta_1 = \frac{1}{4R}$ and $R > 0$ is a large number, with $R \notin \{A_j, B_j\}$. $M_2$ follows the circle $\zeta = \varepsilon[e^{i\theta} - \frac{1}{4}]$ from $\theta = -\frac{\pi}{2}$ to $\theta = \frac{\pi}{2}$, where $\varepsilon$ is a small positive number. (See Figure 2.)

By changing the integral variable, we have

$$(4.2) = \int_L P^v(\zeta; z, t; p) d\zeta,$$

with

$$P^v(\zeta; z, t; p) \ (= P(\zeta; p)) = \frac{p^\zeta}{\Gamma(\zeta+1)} \Phi^v(\zeta; z, t),$$

$$\Phi^v(\zeta; z, t) = \frac{1}{\varphi(G_+(\zeta))} S^v(G_+(\zeta); z, t) G'_+(\zeta).$$

Note that $P(\zeta; p)$ extends analytically with respect to $\zeta$ to the region $\{\text{Re}\zeta > -1\} \setminus (\{A_j\}_{j\in\mathbf{N}} \cup \{B_j\}_{1\leq j\leq J} \cup \{\zeta \leq -\frac{1}{4}\})$ for fixed $z$ and $p$.

LEMMA **4.1** *There are positive numbers $a, t_0, \alpha_0, \kappa_0$ and $\varepsilon_0$ such that, (i) if $\text{Re}\zeta \geq -\frac{1}{4}, |\zeta + \frac{1}{4}| \geq \max\{t_0, a|p|+1\}$ and $\zeta \notin \mathcal{N}$, then*

$$(4.3) \qquad |P(\zeta; p)| \leq C(z, v) r^{-\alpha_0} e^{-\kappa_0 r} |p|^{-\frac{1}{4}},$$

*(ii) if $\text{Re}\zeta = -\frac{1}{4}$ and $|\zeta + \frac{1}{4}| \geq t_0$ or (iii) if $\text{Re}\zeta \geq -\frac{1}{4}$ and $|\zeta + \frac{1}{4}| \leq \varepsilon_0$, then*

$$(4.4) \qquad |P(\zeta; p)| \leq C(z) \left|\zeta + \frac{1}{4}\right|^{-\frac{1}{2m}} \cdot |p|^{-\frac{1}{4}}.$$

Lemma 4.1 (i) implies

$$(4.5) \qquad \int_{M_1} P(\zeta; p) d\zeta \text{ tends to zero as } R \to \infty,$$

whereas Lemma 4.1 (iii) implies

$$(4.6) \qquad \int_{M_2} P(\zeta; p) d\zeta \text{ tends to zero as } \varepsilon \to 0.$$

Figure 2: Integral contours $L$, $M_1$, $M_2$, $N$.

Cauchy's theorem applied to (4.5),(4.6) yields

(4.7) $$\int_L P(\zeta;p)d\zeta = \int_N P(\zeta;p)d\zeta.$$

Here we divide the right-hand side of (4.7) into two parts:

(4.8) $$\int_N P(\zeta;p)d\zeta = S_1(p) + S_2(p),$$

where

$$S_1(p) = i\int_{|t|\leq t_0} P(-\frac{1}{4}+it;p)dt \text{ and } S_2(p) = i\int_{|t|\geq t_0} P(-\frac{1}{4}+it;p)dt.$$

By Lemma 4.1 (iii), we have

$$\begin{aligned} |S_1(p)| &\leq \int_{|t|\leq t_0} |P(-\frac{1}{4}+it;p)|dt \\ &\leq C(z)|p|^{-\frac{1}{4}} \int_{|t|\leq t_0} |t|^{-\frac{1}{2m}}dt \end{aligned}$$

(4.9) $$\leq C(z)|p|^{-\frac{1}{4}}.$$

By Lemma 5.1 (i), we have

$$\begin{aligned} |S_2(p)| &\leq \int_{|t|\geq t_0} |P(-\frac{1}{4}+it;p)|dt \\ &\leq C(z,v)|p|^{-\frac{1}{4}} \int_{|t|\geq t_0} e^{-t}t^{-\frac{5}{4}}dt \end{aligned}$$

(4.10) $$\leq C(z,v)|p|^{-\frac{1}{4}}.$$

Putting (4.8),(4.9),(4.10) together, we have

(4.11) $$|H^v(z,t;p)| = \left|\int_N P(\zeta;p)d\zeta\right| \leq C(z)|p|^{-\frac{1}{4}}.$$

On the other hand, we consider the case where $z$ is in the sector $|\arg z + \frac{\pi}{2}| < \frac{1}{2m-1}\frac{\pi}{2}$. Replacing $+ia_j$ with $-ia_j$, $F_+(\xi)$ with $F_-(\xi)$ and $G_+(\zeta)$ with $G_-(\zeta)$ in the above argument, we can obtain (4.11) in the same fashion.

This completes the proof of Proposition 4.1. □

# 5 Proof of lemmas

In this section we establish the lemmas mentioned previously. We suppose that $z$ is in the sector $|\arg z \pm \frac{\pi}{2}| < \theta_0$, where $\theta_0 := \frac{1}{2m-1}\frac{\pi}{2}$. We set $\zeta = re^{i\theta}$.

## 5.1 Proof of Lemma 3.3

We only give the proof for the case where the integral contour in (3.4) is $\Gamma_+$. We can prove the lemma in the case of $\Gamma_-$ in the same fashion. By Lemma 4.1, we obtain

$$\left|\frac{e^{uv}}{\varphi(v)}\right| \leq C(u)\exp\left\{-c_1 r^{\frac{2m}{2m-1}}\sin\theta_0\right\},$$

and

$$\left|q^{-F_+(v)-1}\right| \leq c|q|^{-\frac{3}{4}}\exp\left\{-c_2 r^{\frac{2m}{2m-1}}[\cos 2\theta_0 \cdot \log|q| - \sin 2\theta_0 \cdot |\arg q|]\right\},$$

for $v \in \Gamma_+$, where $c_1$ is as in Lemma 3.1, $c_2$ is as in Lemma 3.2. Then

$$\left|e^{uv}\frac{q^{-F_+(v)-1}}{\varphi(v)}\right| \leq C(u)|q|^{-1}$$
$$\cdot \exp\left\{-r^{\frac{2m}{2m-1}}\left[c_1\sin\theta_0 + c_2\cos 2\theta_0 \cdot \log|q| - c_2\sin 2\theta_0 \cdot |\arg q|\right]\right\}.$$

Now if we set

$$V = \left\{q \in \mathbf{C}\,;\, |q| \geq \exp\left\{-\frac{1}{3}\frac{\sin\theta_0}{\cos 2\theta_0}\frac{c_1}{c_2}\right\} \text{ and } |\arg q| \leq \min\left\{\frac{\pi}{2}, \frac{1}{3}\frac{\sin\theta_0}{\sin 2\theta_0}\frac{c_1}{c_2}\right\}\right\},$$

then we obtain

$$\left|e^{uv}\frac{q^{-F_+(v)-1}}{\varphi(v)}\right| \leq C(u)e^{-cr^{\frac{2m}{2m-1}}} \quad \text{on } V.$$

Note that $0 < \exp\{-\frac{1}{3}\frac{\sin\theta_0}{\cos 2\theta_0}\frac{c_1}{c_2}\} < 1$. Since the integrand in (3.4) satisfies the above inequality and is a holomorphic function of $q$ on $V$, we can obtain Lemma 3.3. □

## 5.2 Proof of Lemma 4.1

In order to prove Lemma 4.1, we prepare the following two lemmas. We write $P(\zeta;p) = \frac{p^\zeta}{\Gamma(\zeta+1)}\Phi^v(\zeta;z,t)$, where $\Phi^v(\zeta;z,t) = \frac{1}{\varphi(G_\pm(\zeta))}S^v(G_\pm(\zeta);z,t) G'_\pm(\zeta)$.

**LEMMA 5.1** *Let $h$ be a real number. If either* (i) $\operatorname{Re}\zeta = h$ *or* (ii) $\operatorname{Re}\zeta \geq h$ *and $r \geq be|p|$ with $b \geq 1$, then*

(5.1) $$\left|\frac{p^\zeta}{\Gamma(\zeta+1)}\right| \leq c r^{-\frac{1}{2}-h} b^{h-r\cos\theta} |p|^h \exp\{r\sin\theta \cdot \theta\}.$$

**LEMMA 5.2** *There are positive constants $r_0, \beta_0, \varepsilon_0$ such that,* (i) *if* $\operatorname{Re}\zeta \geq -\frac{1}{4}, |\zeta| \geq r_0$ *and $\zeta \notin \mathcal{A}$, then*

(5.2) $$|\Phi^v(\zeta;z,t)| \leq C(z,v) r^{-\beta_0} \exp\left\{-\pi r \frac{\cos(\theta_0 + \frac{\pi}{2} - |\theta|)}{\cos\theta_0}\right\},$$

*and* (ii) *if $\operatorname{Re}\zeta \geq -\frac{1}{4}$ and $|\zeta + \frac{1}{4}| \leq \varepsilon_0$, then*

(5.3) $$|\Phi^v(\zeta;z,t)| \leq C(z)\left|\zeta + \frac{1}{4}\right|^{-\frac{1}{2m}}.$$

By using the above two lemmas, we obtain Lemma 4.1 as follows.

First we consider the case (i). By Lemma 5.1 (ii), where we put $h = -\frac{1}{4}$, and Lemma 5.2 (ii), if $\operatorname{Re}\zeta \geq -\frac{1}{4}, |\zeta + \frac{1}{4}| \geq \max\{r_0 + \frac{1}{4}, be|p|\}$ $(b \geq 1)$ and $\zeta \notin \mathcal{N}$, then we have

$$\begin{aligned}
|P(\zeta;p)| &= \left|\frac{p^\zeta}{\Gamma(\zeta+1)}\right| \cdot |\Phi^v(\zeta;z,t)| \\
&\leq C(z,v) r^{-\beta_0 - \frac{1}{4}} b^{-\frac{1}{4}-r\cos\theta} \\
&\quad \cdot \exp\left\{-\pi r \frac{\cos(\theta_0 + \frac{\pi}{2} - |\theta|)}{\cos\theta_0} + r\sin\theta\cdot\theta\right\} \cdot |p|^{-\frac{1}{4}} \\
&\leq C(z,v) r^{-\beta_0 - \frac{1}{4}} \\
&\quad \cdot \exp\left\{-r\left[\log b\cdot\cos\theta + \pi\frac{\cos(\theta_0 + \frac{\pi}{2} - |\theta|)}{\cos\theta_0} + \sin\theta\cdot\theta\right]\right\} \cdot |p|^{-\frac{1}{4}}.
\end{aligned}$$

We consider the case $\theta \geq 0$. We put $p = \exp\{\frac{\pi}{\cos\theta_0}\}$, then we can obtain

$$\log b \cdot \cos\theta + \pi \frac{\cos(\theta_0 + \frac{\pi}{2} - \theta)}{\cos\theta_0} \geq \frac{1}{\cos\theta_0}\frac{\pi}{2},$$

for $0 \leq \theta \leq \frac{\pi}{2} + \frac{1}{4}(\frac{1}{\cos\theta_0} - 1)$. Therefore if $0 \leq \theta \leq \frac{\pi}{2} + \frac{\pi}{4}(\frac{1}{\cos\theta_0} - 1)$, then we have

$$|P(\zeta;p)| \leq C(z,v)r^{-\beta_0 - \frac{1}{4}} \exp\left\{-r\left[\frac{1}{\cos\theta_0}\frac{\pi}{2} - \sin\theta \cdot \theta\right]\right\}$$

(5.4) $$\leq C(z,v)r^{-\beta_0 - \frac{1}{4}} \exp\left\{-r\left[\frac{\pi}{4}\left(\frac{1}{\cos\theta_0} - 1\right)\right]\right\}.$$

If $-\frac{\pi}{2} - \frac{\pi}{4}(\frac{1}{\cos\theta_0} - 1) \leq \theta \leq 0$, then we obtain the inequality (5.4) in the same fashion. Here if we put $a = \exp\{\frac{\pi}{\cos\theta_0} + 1\}$, $\alpha_0 = \beta_0 + \frac{1}{4}$, $t_0 = \max\{r_0 + \frac{1}{4}, [4\tan\{\frac{1}{4}(\frac{1}{\cos\theta_0} - 1)\}]^{-1}\}$ and $\kappa_0 = \frac{\pi}{4}(\frac{1}{\cos\theta_0} - 1)$, we have (4.3) under the condition (i).

Next we consider the cases (ii), (iii). Since $P(\zeta;p)$ is continuous in $\zeta$ on the set $\{\text{Re}\zeta = -\frac{1}{4} \text{ and } \zeta \neq -\frac{1}{4}\}$, we obtain (4.4) under the condition (ii) by (5.1),(5.3). On the other hand, if we put $\varepsilon_0 = \frac{1}{2}(\min\{A_1, B_1\} + \frac{1}{4})$, then we obtain (4.4) under the condition (iii). In fact $\frac{p^\zeta}{\Gamma(\zeta+1)}$ is continuous in $\zeta$ on the set $\{\text{Re}\zeta \geq -\frac{1}{4}, \text{ and } |\zeta + \frac{1}{4}| \leq \varepsilon_0\}$. □

## 5.3 Proof of Lemma 5.1

When $\text{Re}\zeta \geq h$, we have

$$\frac{1}{\Gamma(\zeta+1)} = \frac{1}{(2\pi)^{\frac{1}{2}}}\frac{e^\zeta}{\zeta^{\zeta+\frac{1}{2}}}\left\{1 + O(\zeta^{-1})\right\}.$$

by Stirling's formula. Hence

$$\left|\frac{p^\zeta}{\Gamma(\zeta+1)}\right| < cr^{-\frac{1}{2}}\exp\{r\cos\theta\log(e|p|r^{-1}) + r\sin\theta \cdot \theta\}$$

$$= cr^{-\frac{1}{2}}b^{-r\cos\theta} \cdot \exp\{r\cos\theta \cdot \log(be|p|r^{-1}) + r\sin\theta \cdot \theta\}.$$

If the condition (i) is satisfied, then

$$\left|\frac{p^\zeta}{\Gamma(\zeta+1)}\right| < cr^{-\frac{1}{2}}b^{-r\cos\theta}\cdot\exp\{h\log(be|p|r^{-1})+r\sin\theta\cdot\theta\}$$
$$< cr^{-\frac{1}{2}-h}b^{h-r\cos\theta}|p|^h\cdot\exp\{r\sin\theta\cdot\theta\}.$$

If the condition (ii) is satisfied, then

$$r\cos\theta\log(be|p|r^{-1}) \leq h\log(be|p|r^{-1}),$$

and the result follows as before. □

## 5.4 Proof of Lemma 5.2

First we consider the case (i). We define the sectors $V_\pm \subset \mathbf{C}\times\mathbf{C}\times\mathbf{R}$ by

$$V_\pm = \left\{(\xi;z,t);\ \left|\arg\xi\mp\frac{\pi}{2}\right|\leq\frac{2m-2}{2m-1}\frac{\pi}{2}\ \text{and}\ \left|\arg z\mp\frac{\pi}{2}\right|<\frac{1}{2m-1}\frac{\pi}{2}\right\}.$$

We need the following lemma about the behavior of $S^v(\cdot;z,t)$ at infinity on $V_\pm$:

LEMMA 5.3 *If $(\xi;z,t)$ is in the sectors $V_\pm$, then there is a non-zero constant $A^v(z)\in\mathbf{C}$ depending on $z$ and $v$ such that*

$$(5.5)\qquad \lim_{|\xi|\to\infty}\xi^{2mv+2m}S^v(\xi;z,t)=A^v(z).$$

By the above lemma, we have

$$(5.6)\qquad |S^v(\xi;z,t)|\leq C(z,v)|\xi|^{-2mv-2m}.$$

By Lemma 3.1, we have the following: If $|\arg\zeta|\leq\frac{\pi}{2}$ and $\zeta\notin\mathcal{N}$, then

$$(5.7)\qquad \left|\frac{1}{\varphi(G_\pm(\zeta))}\right|\leq Cr^{\frac{m-1}{2m}}\exp\left\{-\pi r\frac{\cos(\theta_0+\frac{\pi}{2}-|\theta|)}{\cos\theta_0}\right\}.$$

Since $G'_\pm(\zeta) = \pm\frac{2m-1}{2m}i[\frac{1}{c_2}(\zeta+\frac{1}{4})]^{-\frac{1}{2m}}$, we have

$$(5.8)\qquad |G'_\pm(\zeta)|\leq C\left|\zeta+\frac{1}{4}\right|^{-\frac{1}{2m}}.$$

If we put $\beta_0 = (2m-1)v + (2m + \frac{1}{m} - \frac{3}{2})$, we have (5.2) under the condition (i) by (5.6), (5.7), (5.8). We remark that the above value of $\beta_0$ is best possible.

Next we consider the case (ii). If we put $\varepsilon_0$ as in Subsection 5.2, we have (5.3) by (5.8) under the condition (ii). In fact $[\varphi(G_\pm(\zeta))]^{-1}$ and $S^v(G_\pm(\zeta); z, t)$ are continuous in $\zeta$ on the set $\{\zeta; \operatorname{Re}\zeta \geq -\frac{1}{4}, \text{ and } |\zeta + \frac{1}{4}| \leq \varepsilon_0\}$. □

## 5.5 Proof of Lemma 5.3

First we consider the case $(\xi; z, t) \in V_+$. Writing $\xi = \rho e^{i\alpha}$, with $\rho > 0$ and $|\alpha - \frac{\pi}{2}| \leq \frac{2m-2}{2m-1}\frac{\pi}{2}$, and changing the integral variable, we have

$$\begin{aligned} S^v(\xi; z, t) &= \int_0^\infty e^{it\tau} e^{-x^{2m}\tau} e^{\xi z \tau^{\frac{1}{2m}}} \tau^v d\tau \\ &= \frac{2m}{\rho^{2mv+2m}} \int_0^\infty \exp\{it\frac{u^{2m}}{\rho^{2m}} - x^{2m}\frac{u^{2m}}{\rho^{2m}} + e^{i\alpha}zu\} u^{2mv+2m-1} du. \end{aligned}$$

Since $\operatorname{Re}(e^{i\alpha}z) < 0$ on $V_+$,

$$\int_0^\infty \exp\{it\frac{u^{2m}}{\rho^{2m}} - x^{2m}\frac{u^{2m}}{\rho^{2m}} + e^{i\alpha}zu\} u^{2mv+2m-1} du$$

$$\to \int_0^\infty e^{e^{i\alpha}zu} u^{2mv+2m-1} du = \frac{\Gamma(2mv + 2m)}{[e^{i(\alpha-\pi)}z]^{2mv+2m}} \neq 0$$

as $\rho \to \infty$.

In the case $(\xi; z, t) \in V_-$, we have (5.5) in the same fashion as above. □

## 5.6 Divergence of the formal sum in Remark 3)

If we change the order of the integral and the sum in (1.2) (resp. (1.5)), we obtain the formal sum:

(5.9) $$c \sum_{j=1}^\infty \frac{1}{\varphi'(ia_j)} S_j^v(z, t),$$

where $v = \frac{1}{m}$ (resp. $v = \frac{1}{m} + 1$). Then we have

PROPOSITION **5.1** *The formal sum* (5.9) *is not convergent in the usual sense.*

*Proof.* By Lemma 3.1, we can see that there are $j_0 \in \mathbf{N}$ and $a > 1$ such that
$$\frac{1}{|\varphi'(ia_j)|} \geq a^j$$
for $j \geq j_0$. Moreover by Lemma 5.3, we have
$$|S_j^v(z,t)| \geq cj^{-2mv-2m},$$
on $|\arg z \pm \frac{\pi}{2}| \leq \frac{1}{2m-1}\frac{\pi}{2}$. By the above two inequalities, we see that the formal sum (5.9) is not convergent in the usual sense. □

## 6 Direct proof of Theorem 0.2

By our representation (1.2), we can directly show that the Szegö kernel of $M_m$, off the diagonal, fails to be real analytic and moreover it belongs to $s$-th order Gevrey class $G^s$ for all $s \geq 2m$, but no better. In the case of the Bergman kernel of $D_m$, we can obtain the same result in a similar fashion.

We show that the singularity of the Szegö kernel of $M_m$ is almost equal to that of the function $S_1^{\frac{1}{m}}$, where $S_1^{\frac{1}{m}}$ is the singular solution involving the first eigenfunction for $L_\xi g = 0$ (see §2). Now we suppose that $z$ is in the sector $|\arg z - \frac{\pi}{2}| < \frac{1}{2m-1}\frac{\pi}{2}$. By our representation and the residue formula, we have

(6.1) $$K(z,t) = c \sum_{j=1}^{N-1} \frac{1}{\varphi'(ia_j)} S_j^{\frac{1}{m}}(z,t) + R_N(z,t).$$

Then we have

LEMMA **6.1** *There is a positive constant $C$ independent of $k$ such that*
$$\left|\frac{\partial^k}{\partial t^k} R_N(0+iy, 0)\right| \leq C \frac{\Gamma(2mk+2m+2)}{[|y|(a_N - \varepsilon)]^{2mk+2m+2}},$$
*where $\varepsilon$ is an arbitrary constant with $0 < \varepsilon < a_N - a_{N-1}$.*

By the above lemma and (6.1), we have

$$\frac{\partial^k}{\partial t^k}K(0+iy,0) = c\frac{\partial^k}{\partial t^k}S_1^{\frac{1}{m}}(0+iy,0)\left\{1+O\left(a^{-k}\right)\right\}$$

as $k \to \infty$, where $a > 1$ is a constant. Note that

$$\frac{\partial^k}{\partial t^k}S_j^{\frac{1}{m}}(0+iy,0) = 2mi^k\frac{\Gamma(2mk+2m+2)}{(|y|a_j)^{2mk+2m+2}}$$

(see §2). Since $S_1^{\frac{1}{m}}$ does not belong to $s$-th order Gevrey class $G^s$ for $s < 2m$, we can easily obtain

$$\left|\frac{\partial^k}{\partial t^k}K(0+iy,0)\right| \geq c\frac{\Gamma(2mk+2m+2)}{(|y|a_1)^{2mk+2m+2}}.$$

for sufficiently large $k \in \mathbf{N}$. We can obtain the same result in the case where $|\arg z + \frac{\pi}{2}| < \frac{1}{2m-1}\frac{\pi}{2}$ in the same fashion.

*Proof of Lemma 6.1.* By our representation (1.2), we obtain

$$\frac{\partial^k}{\partial t^k}R_N(0+iy,0) = c(-i)^k 2m\Gamma(2mk+2m+2)|y|^{-2mk-2m-2}I_{k,N},$$

with

$$I_{k,N} = \int_0^\infty e^{-p}\left[\sum_{j=N}^\infty \frac{c_j}{a_j^{2mk+2m+2}}p^{f(j)}\right]dp.$$

Therefore in order to prove Lemma 6.1, it is sufficient to show that

(6.2) $$|I_{k,N}| \leq c(a_N - \varepsilon)^{-2mk},$$

where $\varepsilon$ is an arbitrary constant with $0 < \varepsilon < a_N - a_{N-1}$.

We shall prove the above inequality. By using a similar argument in Sections 3,4 and Fubini's theorem, we have

$$I_{k,N} = \int_C h(v)dv.$$

Here $h(v) = [\varphi(v)]^{-1}v^{-2mk-2m-2}$ and the integral contour $C$ consists of three parts $L_1, L_2, L_3$: $L_1$ follows the half-line from $-\infty$ to $\delta$, $L_2$ follows the circle

Figure 3: Integral contour $C$.

$v = \delta e^{i\theta}$ from $\theta = \pi$ to $\theta = 0$ and $L_3$ follows the half-line from $\delta$ to $\infty$, where $a_{N-1} < \delta < a_N$ and $\delta \notin \{b_j\}$. (See Figure 3.)

Since $h$ is a positive function on $L_3$, we have by Schwarz's inequality

$$\int_{L_3} h(v)dv \leq \left\{\int_\delta^\infty \varphi(v)^{-2} dv\right\}^{\frac{1}{2}} \cdot \left\{\int_\delta^\infty v^{-4mk-4m-4} dv\right\}^{\frac{1}{2}}$$

(6.3)
$$\leq c\delta^{-2mk}.$$

Since $h$ is an even function, we have

(6.4) $$\int_{L_1} h(v)dv \leq \delta^{-2mk}.$$

Moreover we have

$$\left|\int_{L_2} h(v)dv\right| \leq \delta^{-2mk-2m-3} \left|\int_0^\pi \frac{e^{-i(2mk+2m+1)\theta}}{\varphi(\delta e^{i\theta})} d\theta\right|$$

(6.5)
$$\leq c\delta^{-2mk}.$$

Therefore we have

$$\left|\int_C h(v)dv\right| < c\delta^{-2mk},$$

by (6.3),(6.4),(6.5), so we obtain (6.2).

This completes the proof of Lemma 6.1. □

# 7 Singularities on the diagonal

In this section we study the boundary behavior of the Szegö and Bergman kernels of the domain $D_m := \{\text{Im} z_2 > [\text{Re} z_2]^{2m}\}$ ($m = 2, 3, \ldots$) on the diagonal (i.e. $\{(z, w); z = w \in D_m\}$). Here the Szegö kernel $S(z, w)$ of the domain $D \subset \mathbf{C}^n$ is holomorphic in $z$, antiholomorphic in $w$, hermitian symmetric and it agrees with the Szegö kernel on $\partial D \times \partial D$.

The singularities of the Szegö and Bergman kernels at strictly pseudoconvex points are studied in detail (e.g.[15],[4]). We are interested in their behaviors on the neighborhood of the points on $\partial D_m$ at which the Levi form is degenerate. The set of such points is $\{(z_1, z_2) \in \partial D_m; \text{Re} z_1 = 0\}$. On the diagonal the singularities of the Szegö and Bergman kernels are very similar to those of the domain $\{\sum_{j=1}^{n} |z_j|^{2m_j} < 1\} \subset \mathbf{C}^n$ ($m_j \in \mathbf{N}, j = 1, \ldots, n$). The detailed explanation is given in [24]. It is important that the singularities are expressed clearly by using the polar coordinates $(\omega, \rho)$.

PROPOSITION 7.1 *For the domain $D_m = \{\rho(z) := \mathrm{Im}z_2 - [\mathrm{Re}z_1]^{2m} > 0\}$ ($m = 2, 3, \ldots$), the Szegö kernel $S$ and the Bergman kernel $B$ take the following forms:*

$$S(z,z) = \frac{\Phi(\omega(z))}{\rho(z)^{1+\frac{1}{m}}}; \quad B(z,z) = \frac{\Phi^B(\omega(z))}{\rho(z)^{2+\frac{1}{m}}},$$

*where $\omega = [\mathrm{Re}z_1] \cdot [\mathrm{Im}z_2]^{-\frac{1}{2m}}$. Moreover $\Phi(\omega)$ and $\Phi^B(\omega)$ are real analytic functions on $(-1, 1)$ and they have the following estimates:*

$$\frac{c_1}{[1-|\omega|]^{1-\frac{1}{m}}} < \Phi(\omega), \Phi^B(\omega) < \frac{c_2}{[1-|\omega|]^{1-\frac{1}{m}}},$$

*where $c_1, c_2$ are positive constants.*

*Remark.* From the above proposition, $\Phi(\omega), \Phi^B(\omega)$ are bounded on the region $\{\mathrm{Im}z_2 > \alpha[\mathrm{Re}z_1]^{2m}\}$ ($\alpha > 1$), but they are not elsewhere.

*Proof.* We write $A \approx B$ to indicate that there are positive constants $c_1, c_2$ such that $c_1 A < B < c_2 A$.

By [30],[32], the Szegö kernel $S$ and the Bergman kernel $B$ of the domain $D_m$ can be expressed in the following:

$$S(z,z) = c \int_0^\infty e^{-b\tau} P(a\tau^{\frac{1}{2m}}) \tau^{\frac{1}{m}} d\tau;$$

$$B(z,z) = c \int_0^\infty e^{-b\tau} P(a\tau^{\frac{1}{2m}}) \tau^{1+\frac{1}{m}} d\tau,$$

where $a = \mathrm{Re}z_1$, $b = \mathrm{Im}z_2$ and $P$ is as in (3.2) in Section 3. Now we introduce the new variable $\omega = a \cdot b^{-\frac{1}{2m}}$ in the above equations. Then we have

$$S(z,z) = \frac{\Phi(\omega(z))}{\rho(z)^{1+\frac{1}{m}}}; \quad B(z) = \frac{\Phi^B(\omega(z))}{\rho(z)^{2+\frac{1}{m}}},$$

with

$$\Phi(\omega) = c[1-\omega^{2m}]^{1+\frac{1}{m}} \int_0^\infty e^{-s} P(\omega s^{\frac{1}{2m}}) s^{\frac{1}{m}} ds;$$

$$\Phi^B(\omega) = c[1-\omega^{2m}]^{2+\frac{1}{m}}\int_0^\infty e^{-s}P(\omega s^{\frac{1}{2m}})s^{1+\frac{1}{m}}ds.$$

In order to estimate the function $P$, we use a part of Haslinger's results in [20].

LEMMA **7.1** ([20]) *For $u \geq 0$,*

$$P(u) \approx [1+u^{2m-2}]e^{u^{2m}}.$$

From the above lemma, we have

$$\begin{aligned}
\Phi(\omega) &\approx [1-\omega^{2m}]^{1+\frac{1}{m}}\int_0^\infty e^{-[1-\omega^{2m}]s}(1+\omega^{2m-2}s^{1-\frac{1}{m}})s^{\frac{1}{m}}ds \\
&= \Gamma(1+\frac{1}{m}) + \frac{\omega^{2m-2}}{[1-\omega^{2m}]^{1-\frac{1}{m}}}
\end{aligned}$$

(7.1) $$\approx \frac{1}{[1-|\omega|]^{1-\frac{1}{m}}}$$

In the same fashion, we have

(7.2) $$\Phi^B(\omega) \approx \frac{1}{[1-|\omega|]^{1-\frac{1}{m}}}.$$

This completes the proof of Proposition 7.1. □

## ACKNOWLEDGMENTS

I would like to express my deepest gratitude to my advisors, Prof. Hiroki Tanabe and Prof. Kazuo Okamoto, for their constant encouragement. I am also indebted to Prof. Katsunori Iwasaki who carefully read the manuscript and supplied many corrections.